\input amstex
\input amsppt.sty   
\hsize 30pc
\vsize 47pc
\magnification=\magstep1
\def\nmb#1#2{#2}         
\def\idx{}               
\def\ign#1{}             

\redefine\o{\circ}
\define\X{\frak X}

\define\de{\delta}

\define\si{\sigma}

\define\ph{\varphi}

\define\La{\Lambda}
\define\Si{\Sigma}

\define\row#1#2#3{#1_{#2},\ldots,#1_{#3}}
\redefine\D{{\Cal D}}

\define\ev{\operatorname{ev}}
\define\x{\times}
\define\Fl{\operatorname{Fl}}

\def\today{\ifcase\month\or
 January\or February\or March\or April\or May\or June\or
 July\or August\or September\or October\or November\or December\fi
 \space\number\day, \number\year}
\topmatter
\title Radon transform and curvature
\endtitle
\author  Peter W. Michor  \endauthor
\affil
Erwin Schr\"odinger Institute for Mathematical Physics, 
and Universit\"at Wien
\endaffil
\address
Erwin Schr\"odinger Institute for Mathematical Physics,
Pasteurgasse 4/7, A-1090 Wien, Austria.
\endaddress
\address
Institut f\"ur Mathematik, Universit\"at,
Strudlhofgasse 4, A-1090 Wien, Austria
\endaddress
\email michor\@pap.univie.ac.at \endemail
\abstract We interpret the setting for a Radon transform as a 
submanifold of the space of generalized functions, and compute its 
extrinsic curvature: it is the Hessian composed with the Radon 
transform.
\endabstract
\endtopmatter

\document

\subhead\nmb0{1}. The general setting \endsubhead
Let $M$ and $\Si$ be smooth finite dimensional manifolds. Let 
$m=\dim(M)$.
A linear mapping $R:C^\infty_c(M)\to C^\infty(\Si)$ is called a 
(generalized) Radon transform if it is given in the following way:
To each point $y\in \Si$ there corresponds a submanifold $\Si_y$ of 
$M$ and a density $\mu_y$ on $\Si_y$, and the operator $R$ is given 
by 
$$R(f)(y):= \int_{\Si_y}f(x)\mu_y(x).$$
We will express this situation in the following way.

Let $\D(M):=C^\infty_c(M)$ be the space of smooth functions with 
compact support on $M$, and let $\D'(M)=C^\infty_c(M)'$ be the 
locally convex dual space. Note that the space $C^\infty(|\La^m|(M))$ 
of smooth densities on $M$ is canonically contained and dense in 
$\D'(M)$.

Now suppose that we are given a smooth mapping $\si:\Si\to \D'(M)$. 
By the smooth uniform boundedness principle
(see \cite{Fr\"olicher, Kriegl, p\. 73} or \cite{Kriegl, Michor, 4.11})
the mapping $\si:\Si\to L(\D(M),\Bbb R)$ is smooth if and only 
if the composition with the evaluation 
$\ev_f:L(\D(M),\Bbb R)\to \Bbb R$ is smooth for each $f\in \D(M)$, 
i\.e\. $R_{\si}(f):\Si\to\Bbb R$ is smooth for each $f$.
Then we have an associated \idx{\it Radon transform} given by 
$$R_{\si}(f)(y):= \langle  \si(y), f\rangle.$$
Clearly the Radon transform $R_\si:C^\infty_c(M)\to C^\infty(\Si)$ is 
injective if and only if the subset $\si(\Si)\subset\D'(M)$ separates 
points on $C^\infty_c(M)$, and the kernel of $R_\si$ is the 
annihilator of $\si(\Si)$ in $C^\infty_c(M)$.
We will assume in the sequel that $\si:\Si\to \D'(M)$ is an embedding 
of a smooth finite dimensional embedded submanifold of the locally 
convex vector space $\D'(M)$, but the Radon transform itself is 
defined also in the more general setting of a smooth mapping.

All examples of Radon transforms mentioned in these proceedings fit 
into the setting explained above. A trivial example is the Dirac 
embedding $\de:M\to \D'(M)$ associating to each point $x\in M$ the 
Dirac measure $\de_x$ at that point. It's associated Radon transform 
is the identity for functions on  $M$, but it's curvature (see below) 
is quite interesting.

\subhead\nmb0{2}. Curvature \endsubhead
We now give the definition of the \idx{\it second fundamental form} or 
the \idx{\it extrinsic curvature} of a finite dimensional submanifold 
$\Si$ of the locally convex space $\D'(M)$. Since we do not want to 
assume the existence of an inner product on (a certain subspace of) 
$\D'(M)$ we consider the \idx{\it normal bundle} $N(\Si):= 
(T\D'(M)|\Si)/T\Si$ and the canonical projection 
$\pi:T\D'(M)|\Si\to N(\Si)$ of vector bundles over $\Si$. The linear 
structure of $\D'(M)$ gives us the obvious flat covariant 
derivative $\nabla_XY$ of two vector fields $X,Y$ on $\D'(M)$, which 
is defined by $(\nabla_XY)(\ph)=dY(\ph).X(\ph)$.
For (local) vector fields $X,Y\in\X(\D'(M))$ on $\D'(M)$ which along $\Si$ 
are tangent to $\Si$ we consider the section $S(X,Y)$ of $N(\Si)$ 
which is given by $S(X,Y)=\pi(\nabla_XY)$. This section depends only 
on $X|\Si$ and $Y|\Si$, since we may consider the flow $\Fl^{X|\Si}_t$ of 
the vector field $X|\Si$ on the finite dimensional manifold $\Si$ and 
we have $(\nabla_XY)|\Si=\tfrac d{dt}|_{t=0}Y\o \Fl^{X|\Si}_t$. Here 
we consider just the smooth mapping $Y:\D'(M)\to\D'(M)$. Obviously 
$S(X,Y)$ is $C^\infty(M)$-linear in $X$, and it is symmetric since 
$S(X,Y)-S(Y,X)=\pi(dY.X-dX.Y)=\pi([X,Y])=0$.
So the \idx{\it second fundamental form} or 
the \idx{\it extrinsic curvature} of the submanifold $\Si$ of 
$\D'(M)$ is given by 
$$\gather
S:T\Si\x_\Si T\Si \to N(\Si).\\
S(X,Y)=\pi(\nabla_XY)\text{ for }X,Y\in\X(\Si).
\endgather$$
For $y\in\Si$ the convenient vector space $N_y(\Si)=\D'(M)/T_y\Si$ is 
the dual space of the closed linear subspace 
$\{f\in\D(M):\langle T_y\si.X,f\rangle=0\text{ for all }X\in T_y\Si\}$.

\proclaim{\nmb0{3}. Theorem}
Let $\si:\Si\to \D'(M)$ be a smooth embedding of a finite dimensional 
smooth manifold $\Si$ into the space of distributions on a manifold 
$M$, and let $R_\si:C^\infty_c(M)\to C^\infty(\Si)$ be the associated 
Radon transform. Then the extrinsic curvature of $\si(\Si)$ in 
$\D'(M)$ is the Hessian composed with the Radon transform in the 
sense explained in the proof.
\endproclaim

\demo{Proof}
Since $\si(\Si)$ is an embedded submanifold of finite dimension in 
$\D'(M)$, it is also splitting, and thus for each vector field 
$X\in \X(\Si)$ there exists a (local) smooth extension 
$\tilde X\in\X(\D'(M))$. 
It is not known whether $\D'(M)$ admits smooth 
partitions of unity. The space $C^\infty_c(M)$ of test functions 
admits smooth partitions of unity, see \cite{Kriegl, Michor}.
So we have $T\si\o X = \tilde X\o \si$.

For $y\in\Si$ the normal space $N_y(\Si)= \D'(M)/T_y\si(T_y\Si)$ is 
the dual space of the annihilator of $T_y\si(T_y\Si)$ in 
$C^\infty_c(M)$. A test function $f\in C^\infty_c(M)$ is in this 
annihilator if and only if $\langle T_y\si.X,f\rangle=0$ for all 
$X\in T_y\Si$. Let us choose a smooth curve $c:\Bbb R\to \Si$ with 
$c(0)=y$ and $c'(0)=X$. Then we have
$$\align
\langle T_y\si.X,f\rangle &= \langle \tfrac d{dt}|_0\si(c(t)),f\rangle
     = \tfrac d{dt}|_0 \langle \si(c(t)),f\rangle \\
&= \tfrac d{dt}|_0 R_\si f(c(t)) = d(R_\si f)_y(X).
\endalign$$
So we have $N_y(\Si) = \{f\in C^\infty_c(M): d(R_\si f)_y=0\}'$.

Now we will compute the extrinsic curvature. Let $X,Y\in \X(\Si)$ be 
vector fields, let $\tilde X, \tilde Y$ be smooth extensions to 
$\D'(M)$, let $y\in \Si$, and choose $f\in C^\infty_c(M)$ with 
$d(R_\si f)_y=0$. Then we have
$$\align
\langle S(X,Y)(y),f\rangle 
     &= \langle (\nabla_{\tilde X}\tilde Y)(\si(y)),f\rangle \\
&= \langle d\tilde Y(\si(y)).\tilde X(\si(y)),f\rangle \\
&= \langle d\tilde Y(\si(y)).d\si(y).X(y),f\rangle \\
&= \langle d(\tilde Y\o \si)(y).X(y),f\rangle \\
&= \langle d(d\si.Y)(y).X(y),f\rangle, \\
Y(R_\si f) &= d(R_\si f).Y = \tfrac d{dt}|_0 R_\si f\o \Fl^Y_t \\
&= \tfrac d{dt}|_0 \langle \si\o\Fl^Y_t,f\rangle 
     = \langle d\si.Y,f\rangle, \\
XY(R_\si f)(y) &= \tfrac d{dt}|_0 (Y(R_\si f))(\Fl^X_t(y)) 
     = \tfrac d{dt}|_0 \langle (d\si.Y)(\Fl^X_t(y)), f\rangle \\
&= \langle d(d\si.Y).X(y),f\rangle = \langle S(X,Y)(y),f\rangle. 
\endalign$$
So $\langle S(X,Y)(y),f\rangle$ is the Hessian of $R_\si f$ at $y$ 
applied to $(X(y),Y(y))$.
\qed\enddemo

\enddocument